\magnification=\magstep1 
\baselineskip=12pt 

\def\ind {{\lim \limits_{\longrightarrow}}} 
\def \pro {{\lim\limits_{\longleftarrow}}} 
\input amstex 
\documentstyle{amsppt} 
\overfullrule=0pt 
\vsize=500pt 

\document

\def\a{{\alpha}}

\def\a{{\alpha}}

\def\eps{{\varepsilon}}

\def\1b{{\bold 1}}

\def\Cb{{\bold C}}

\def\Spec{\roman{Spec\,}} 
\def\Spf{\roman{Spf\,}}

\def\Db{{\bold{D}}} 
\def\Zb{{\bold{Z}}}

\def\ch{{\text{ch}}} 
 
\def\Algb{{\bold{Alg}}} 
\def\Affb{{\bold{Aff}}} 
\def\Schb{{\bold{Sch}}} 
\def\Shb{{\bold{Sh}}}

\def\Grc{{\Gc r}} 
\def\Ischb{{\bold{Isch}}}

\def\Lrsb{{\bold{Lrs}}} 

\def\Setsb{{\bold{Sets}}} 
\def\Affb{{\bold{Aff}}} 
\def\Ischb{{\bold{Isch}}}

\def\det{\text{det}} 
\def\max{{\text{max}}} 
\def\Hom{\roman{Hom}} 
\def\Det{\Dc et} 
\def\Zar{\roman{Zar}}

\def\Id{\roman{Id}}

\def\Res{\text{Res}}

\def\max{\text{max}\,}

\def\Pic{\text{Pic}}

\def\GL{{\roman{GL}}} 
\def\SL{{\roman{SL}}} 
 
\def\Mor{\text{Mor}} 
\def\Map{\text{Map}}

\def\AA{{\Bbb A}} 
 
\def\CC{{\Bbb C}}

\def\HH{{\Bbb H}}

\def\ZZ{{\Bbb Z}}

\def\Ac{{\Cal A}} 
\def\Bc{{\Cal B}} 
 
\def\Dc{{\Cal D}} 
\def\Ec{{\Cal E}} 
\def\Fc{{\Cal F}} 
\def\Gc{{\Cal G}}

\def\Lc{{\Cal L}} 
\def\Mc{{\Cal M}} 
\def\Nc{{\Cal N}} 
\def\Oc{{\Cal O}}

\def\Rc{{\Cal R}} 
\def\Sc{{\Cal S}} 
\def\Tc{{\Cal T}}

\def\Xc{{\Cal X}}

\def\and{{\quad\text{and}\quad}}

\def\ts{\textstyle}

\def\qed{\hfill $\sqcap \hskip-6.5pt \sqcup$}

\def\Wedge{{\ts\bigwedge}}

\def\u1{{\underline 1}}

\centerline{\bf FORMAL LOOPS II: A LOCAL RIEMANN-ROCH THEOREM} 
\centerline{\bf FOR DETERMINANTAL GERBES} 

\vskip .5cm 

\centerline{\bf M. Kapranov, E. Vasserot}

\vskip 1cm

\noindent {\bf (0.1)} The goal of this paper  and the next one [KV2] 
is to relate three subjects of recent interest: 

\vskip .3cm 

\noindent { (A)} The theory of sheaves of chiral differential 
operators (CDO), see [GMS1-2]. A sheaf of CDO on a complex manifold $X$ 
is a sheaf of graded vertex algebras with certain conditions on the 
graded components. As shown in {\it loc. cit.}, locally on $X$ such an object 
always exists and is unique up to an isomorphism but the isomorphism 
not being canonical, the global situation is similar to the behavior of 
spinor bundles on a Riemannian manifold. This is expressed by saying 
that sheaves of CDO form a gerbe $\Cal{CDO}_X$. A global object exists 
if and only if the characteristic class 
$$ ch_2(X) ={1\over 2} c_1^2(X)-c_2(X) \leqno (0.1.1)$$ 
vanishes. Manifolds with this property are known as 
$MU\langle 8 \rangle$-manifolds in homotopy theory. 

\vskip .2cm 

\noindent { (B)} The theory of the group 
$GL(\infty)$ developed by Sato and others [PS], in particular, of the 
determinantal central extension of $GL(\infty)$ and of the semi-infinite 
Grassmann manifold (the Sato Grassmannian). In  fact, what appears here 
is a nonlinear version of this theory. 

\vskip .2cm 

\noindent {(C)} The refinement of the Grothendieck-Riemann-Roch theorem 
found by Deligne [De1]. Recall that if 
$\pi: {\Cal X}\to B$ is a smooth projective morphism of
 algebraic varieties 
and ${\Cal E}$ is an algebraic
 vector bundle on $\Cal X$,  then Grothendieck-Riemann-Roch 
describes the Chern classes $c_i(R\pi_*({\Cal E}))$ in terms of the 
cohomological direct images of the Chern classes of ${\Cal E}$ and of 
the relative tangent bundle $\Theta_{{\Cal X}/B}$. When $i=1$, 
the first Chern class in question can be seen as the class of the 
determinantal line bundle $\det\, R\pi_*({\Cal E})$ in the Picard group of 
$B$ and Deligne posed the problem of identifying this bundle itself 
up to a canonical isomorphism. He solved this problem for the 
case when $\pi: {\Cal X}\to B$ is a family of curves. Our theory 
uncovers a new, deeper level of this subject: that of a 
``local'' version of Deligne-Riemann-Roch describing determinantal 
gerbes, not determinantal line bundles. 

\vskip .3cm 

\noindent {\bf (0.2)} The relation (A)-(B) proceeds via the ind-scheme 
${\Cal L}(X)$ of formal loops introduced in [KV1]. It contains the 
scheme ${\Cal L}^0(X)$ of formal arcs, see [DL]. They are 
algebro-geometric analogs of the Frechet manifolds 
$$LX = C^\infty(S^1, X),  \quad L^0X = \roman{ Hol}(D^2, X), \leqno (0.2.1)$$ 
where $D^2$ is the unit disk $\{|z|\leq 1\}$ and Hol stands for the space 
of holomorphic maps smooth on the boundary. 

\vskip .2cm 

Sato's theory of $GL(\infty)$ can be developed in two versions. 
The formal version works with a linearly compact topological 
vector space $V$ such as the space $\CC((t))$ with the $t$-adic topology. 
The group $GL(\infty)$ is then incarnated as the group of 
continuous automorphisms of $V$. This is the version we adopt in the main body 
of the paper. The Hilbert version starts with a 
Hilbert space $H$ equipped with a polarization: a class of splittings 
$H= H_+\oplus H_-$ defined up to some equivalence, see [PS]. The group 
$GL(\infty)$ is incarnated as the group of bounded linear automorphisms 
preserving the polarization. In both cases one has a Grassmann-type variety 
$\Cal G$ (the Sato Grassmannian) 
and a determinantal line bundle $\Delta$ on ${\Cal G}\times 
{\Cal G}$ making $\Cal G$ into a set of objects of a $\CC^\times$-gerbe acted 
 upon by $GL(\infty)$. 

\vskip .2cm 

A nonlinear version of the theory should involve infinite-dimensional 
manifolds with a $GL(\infty)$-structure in the tangent bundle. 
In the sequel to this  paper [KV2] we develop a convenient formalization 
of this idea in the algebro-geometric setting. 
The corresponding objects are locally compact smooth ind-schemes, 
introduced in \cite{KV1}. Their tangent spaces possess
 a $GL(\infty)$-structure 
in the formal version above. In particular, for such an ind-scheme $Y$ 
we have the relative Sato Grassmannian ${\Cal G}\to Y$ 
and a determinantal line bundle $\Delta$ 
on ${\Cal G}\times_Y {\Cal G}$ which gives rise to an 
 ${\Cal O}^\times_Y$-gerbe 
$\Det_Y$. The ``determinantal anomaly'' of $Y$ is the class of 
$\Det_Y$ in the group $H^2(Y, {\Cal O}_Y^\times)$ classifying 
${\Cal O}_Y^\times$-gerbes. An example is provided by $Y={\Cal L}X$. 
It turns out that it is precisely this gerbe that governs 
sheaves of CDO. In other words, the anomaly in constructing CDO 
is precisely the determinantal anomaly for the loop space [KV2].

Our main result, Theorem 5.3.1, has the following consequence. The class 
$[\Det_{{\Cal L}X}]\in H^2({\Cal L}X, {\Cal O}^\times)$ is equal to 
 the image of the characteristic class (0.1.1) under the transgression map 
(tame symbol) 
$$\tau: H^2(X, K_2({\Cal O}_X)) \to H^2({\Cal L}X, {\Cal O}^\times).$$ 

\vskip .3cm 

\noindent {\bf (0.3)} The above identification of $[\Det_{{\Cal L}X}]$ 
can be seen as a particular case of a Riemann-Roch-type result for 
determinantal gerbes, and it is this Riemann-Roch theorem for gerbes 
that constitutes the main 
result of this paper.  It is easier to explain the situation in the 
$C^\infty$-framework (0.2.1). Let $\Sigma$ be an oriented $C^\infty$-manifold 
diffeomorphic to $S^1$ and $E$ a $C^\infty$ complex vector bundle on 
$\Sigma$. We have then the space of smooth sections $\Gamma^{C^\infty} 
(\Sigma, E)$ and its Hilbert space completion $H=\Gamma^{L^2}(\Sigma, E)$. 
The latter can be defined using any Riemannian metric on $\Sigma$ and 
Hermitian metric on $E$ and is independent, as a topological 
vector space, of the choices. It has a canonical 
$GL(\infty)$-structure. To define it, one realizes $\Sigma$ as the boundary 
of a holomorphic disk $D$ and extends $E$ to a holomorphic bundle $\Cal E$ 
on $D$. Then one takes $H_+= \Gamma^{\roman hol}(D, {\Cal E})$ (holomorphic 
sections) and the resulting $GL(\infty)$-structure is independent on the 
choice of $D,{\Cal E}$. Now, let $p: \Sigma\to B$ is a smooth $S^1$-fibration 
over a $C^\infty$ (or Frechet) manifold $B$ and suppose that fibers of $p$ are 
oriented. If $E$ is a smooth 
vector bundle on $\Sigma$, we have a ``bundle of vector spaces with 
$GL(\infty)$-structure'' on $B$, namely the $L^2$-direct image 
$p_*^{L^2}(E)$. It gives a gerbe $\Det(p_*E)$ 
with band ${\Cal O}^\times_B$, the sheaf of invertible $C^\infty$-functions 
and hence gives a class 
$$[\Det(p_*E)] \in H^2(B, {\Cal O}_B^\times). \leqno (0.3.1)$$ 
Denoting by $\delta: H^2(B, {\Cal O}_B^\times)\to H^3(B, {\ZZ})$ 
the coboundary map of the exponential sequence, we get a class 
$\delta[\Det(p_*E)]\in H^3(B, {\ZZ})$. For $B=LX$ the tangent bundle 
$\Theta_{LX}$ can be represented as $p_* {\roman ev}^*(\Theta_X)$ where 
in the diagram 
$$LX\buildrel p\over\longleftarrow S^1\times LX 
\buildrel{\roman ev}\over\longrightarrow X \leqno (0.3.2)$$ 
$p$ is the projection and ev is the evaluation map. The class 
$[\Det_{LX}]\in H^2(LX, {\Cal O}^\times)$ is a particular case of (0.3.1). 

\vskip .2cm 

Now, the ``local Riemann-Roch'' for a smooth circle fibration 
$p: \Sigma\to B$ says that after inverting 2, one has 
$$\delta[\Det(p_*E)] = \int_{\Sigma/B}\biggl({1\over 2} c_1^2(E)-c_2(E)
\biggr)\in 
H^3(B, {\ZZ})\otimes {\ZZ}\bigl[{1\over 2}\bigr], \leqno (0.3.3)$$ 
where $\int_{\Sigma/B}: H^4(\Sigma, {\ZZ})\to H^3(B, {\ZZ})$ 
is the cohomological pushforward. This theorem 
 and its higher-dimensional version (see below) will be 
proved in a subsequent paper. The  main result of the present paper 
(Theorem 5.3.1) 
is a formal version of (0.3.3) (with formal Laurent series replacing 
functions on $S^1$) for the trivial fibration. 

\vskip .3cm 

\noindent {\bf (0.4)} Theorem 5.3.1 can in fact be seen as a statement 
comparing two central extensions of the loop group $GL_N((t))$ 
(or $L(GL_N)$ in the analytic setting). If $G$ is a simple, simply 
connected algebraic group over $\CC$, then the group of central extensions 
of $G((t))$ by $\CC^*$ is generated by the class of the so-called 
``primitive'' extension which goes back to the works 
of Moore [Mo] and Matsumoto [Ma]. For $G=SL_N$ 
the primitive extension (or the universal central extension with
center $K_2$) 
was related with the second Chern class by S. Bloch [Bl]. It is also well 
known in the theory of Kac-Moody groups [PS] [Ku] that the primitive 
extension for $SL_N((t))$ is precisely the determinantal one. 
Combined together, these observations can be used to prove a 
particular case of (0.3.3): that for $E$ with $\det(E)$ trivial, we have 
$$\delta[\Dc et(p_*(E)] = -\int_{\Sigma/B} c_2(E). \leqno (0.4.1)$$ 
The case of more general reductive groups $G$ is more complicated 
and was studied in [BD] and other papers. Our result identifies the 
determinantal central extension of $GL_N((t))$ with the 
extension coming from $ch_2$. It does not directly follow 
from the cited papers since we need an identification 
as a group ind-scheme, not as a discrete group. 

\vskip .3cm

The relation of (0.3.3) to Deligne-RR can be understood as follows. 
Suppose that $B$ is a complex algebraic variety and $\Sigma$ is 
embedded into a family $\pi: {\Cal X}\to B$ of algebraic curves, 
as in (0.1)(C) so that $\Xc$ is split into the ``inside'' and 
``outside'' parts: $\Xc = \Xc_+ \cup \Xc_-$. Then every 
extension $\Ec_+$ of $E$ to a holomorphic bundle on 
$\Xc_+$ gives an object $[\Ec_+]$ of $\Dc et(p_*E)$. 
Given extensions $\Ec_+, \Ec_-$ to $\Xc_+, \Xc_-$ we have a 
holomorphic bundle $\Ec$ to which Deligne-RR applies. Fixing $\Ec_-$ 
and taking two choices $\Ec'_+, \Ec''_+$, we get two 
holomorphic bundles $\Ec', \Ec''$ and 
$${\roman Hom} _{\Dc et(p_*E)}([\Ec'_+], [\Ec''_+]) = 
{\det(R\pi_* \Ec'')\over \det (R\pi_*\Ec')}. \leqno (0.4.2)$$ 
So we recognize in (0.3.3) the 
beginning of the $ch.Td$ quantity of the Riemann-Roch. 

\vskip .3cm

Put differently, the relation of  the Riemann-Roch theorem to (0.3.3) 
is similar to the relation of the self-duality of the Jacobian of 
a curve to the Cartier self-duality of the ind-group scheme 
$GL_1((t))$ established by C. Contou-Carrere [CC]. In fact, the 
Contou-Carrere symbol plays an important role in our approach. 

\vskip .3cm 

\noindent {\bf (0.5)} If now $p: \Sigma\to B$ 
is a $C^\infty$-fibration of compact oriented $C^\infty$-manifolds 
 of dimension 
$d$ and $E$ is a vector bundle on $\Sigma$, then it is natural to expect 
 that the fibers of 
$p_*^{L^2}(E)$ have some kind of ``$d$-fold polarization'' and give rise to 
a determinantal $d$-gerbe $\Det(p_*E)$ with band ${\Cal O}_B^\times$. 
Although one does not know how to define this $d$-gerbe, one can define 
the class 
$C_1(p_*E) \in H^{d+2}(B, \CC)$ 
which should be the image of $[\Det(p_*E)] \in H^{d+1}(B, {\Cal O}^\times)$ 
 under the composite map 
$$H^{d+1}(B, {\Cal O}^\times) \buildrel \delta_{\exp}\over\longrightarrow 
H^{d+2}(B, \ZZ)\buildrel \otimes \CC\over\longrightarrow 
H^{d+2}(B, \CC).$$ 
This can be done using the Chern-Weil approach and we have 
$$C_1(p_*E) = \int_{\Sigma/B} \biggl[ ch(E) 
 \cdot {\roman Td}(\Theta_{\Sigma/B}\otimes \CC)\biggr]_{d+1}. 
\leqno (0.4.1)$$ 
Note that this statement involves only real manifolds and 
utilizes the component of $ch \cdot{\roman Td}$ quite different from the 
standard GRR. It will be proved in another  paper [BKTV]  devoted to the 
$C^\infty$ rather than formal theory. 

\vskip .3cm

The research of M.K. was partially supported by an NSF grant. 

\vskip 2 cm 

\centerline{\bf 1. Reminder on Grothendieck sites and  gerbes.} 

\vskip .3cm

Recall [Mi] that a Grothendieck site is a category $\Sc$ with fiber products 
and disjoint unions equipped with a family of coverings of objects so that one 
can speak about sheaves on $\Sc$ as functors satisfying descent for coverings. 
Objects $U$ of $\Sc$ will be referred to as open sets and we will write 
$U\subset \Sc$ 
to emphasize the analogy with topological spaces.

Lets $\Sc$, $\Tc$ be two Grothendieck sites. 
 A morphism of sites 
$f:\Sc\to\Tc$ is a functor from the underlying category of $\Tc$ 
to the underlying category of $\Sc$ commuting with disjoint unions, fiber 
products, and taking coverings to coverings. 
If $\Ac$ is a sheaf of rings on $\Tc$, a morphism of ringed sites 
$$f:(\Sc,\Ac)\to(\Tc,\Bc)$$ 
is a pair consisting of a morphism of sites 
$$f_\sharp:\Sc\to\Tc$$ 
and a morphism of sheaves of rings on $\Sc$ 
$$f^\flat:f_\sharp^{-1}(\Bc)\to\Ac.$$ 
If $\Mc$ is a sheaf of $\Bc$-modules on $\Tc$, we denote by 
$$f^\circ \Mc = f_\sharp^{-1} \Mc \otimes _{f_\sharp^{-1} \Bc} \Ac$$ 
its inverse image in the sense of moprhisms  of ringed sites. 
Note that we have a morphism 
$f^\circ  : H^i(\Tc, \Mc) \to H^i(\Sc, f^\circ \Mc)$.

If $\Sc$ is a Grothendieck site 
and $\Cal F$ is a sheaf of abelian groups on $\Sc$, then we can speak of 
$\Cal F$-gerbes (= gerbes with band $\Cal F$, see [Br]). Recall that such a 
gerbe $\Cal G$ 
consists of the following data: 

\itemitem{(1)} A category ${\Cal G}(U)$ given for all open $U\subset \Sc$, 
the restriction functors $r_{UV}: {\Cal G}(U)\to {\Cal G}(V)$ given for any 
morphism 
$V\to U$ 
and natural isomorphisms of functors $s_{UVW}: r_{VW}\circ r_{UV}\Rightarrow 
r_{UW}$ 
given for each $W\to V\to U$ and satisfying the transitivity conditions.

\itemitem{(2)} The structure of ${\Cal F}(U)$-torsor on each 
$\Hom_{{\Cal G}(U)} 
(x,y)$ compatible with the $r_{UV}$ and such that the composition of morphisms 
is bi-additive. 

These data have to satisfy the local uniqueness and gluing properties for which 
we refer to [Br]. 

By a sheaf of $\Cal F$-groupoids we will mean a sheaf of categories $\Cal C$ 
on $\Sc$ 
(so both $\roman{ Ob}\, {\Cal C}$ and 
$\Mor \, {\Cal C}$ are sheaves of sets) 
in which each $\Hom_{{\Cal C}(U)}(x,y)$ is either empty or is made 
into an ${\Cal F}(U)$-torsor 
so that the composition is biadditive. 
A sheaf $\Cal C$  of $\Cal F$-groupoids is called 
locally connected if locally on $S$ all the $\roman{Ob} \, {\Cal C}(U)$ 
and  $\Hom_{{\Cal C}(U)}(x,y)$ 
are nonempty.

Each sheaf of $\Cal F$-groupoids can be seen as a fibered category 
over $\Sc$, in fact it is a pre-stack, see, e.g., [LM]. 
Recall (see, e.g., {\it loc. cit.} Lemma 2.2) that for any pre-stack 
$\Cal C$ there is an associated stack ${\Cal C}^{^\sim}$. 
If $\Cal C$ is a locally connected sheaf of $\Cal F$-groupoids, 
then ${\Cal C}^{^\sim}$ is an $\Cal F$-gerbe. We will refer to 
this procedure as ``gerbification" (analog of sheafification). 

As well known (see, e.g., [Br]), the set formed by 
${\Cal F}$-gerbes up to equivalence is identified with $H^2(\Sc, {\Cal F})$. 
Given an $\Cal F$-gerbe $\Cal G$, we denote by $[{\Cal F}]\in H^2(S, {\Cal F})$ 
its class. Given a sheaf $\Cal C$ of $\Cal F$-groupoids, we denote 
by $[{\Cal C}]$ the class of the corresponding gerbe.

\vskip 2cm

\centerline{\bf 2. Ind-schemes and their cohomology.} 

\vskip .3cm

Let $F$ be a field. By 
$\Schb_F $ we denote the category of $F$-schemes. 
As in [KV1],  by an ind-scheme over $F$  we will mean an ind-object of
 $\Schb_F$ of the form $``\ind"Y^\alpha$ where $\a$ runs in a filtered 
poset $A$, $Y^\alpha$ are quasicompact schemes (not 
necessarily of finite type) 
and the structural maps $Y^\alpha\to Y^{\alpha'}$ are closed 
embeddings. We will say that $Y$ is of countable type if 
it admits a presentation as above with a countable poset $A$. 
Sometimes we use the notation $Y^\infty$ instead of $Y$ to empashize that 
we are dealing with an ind-scheme. We denote by $\Ischb_F$ the category of 
ind-schemes over $F$. In the case $F=\CC$ we write simply
 $\Schb$ and $\Ischb$.

As with any ind-object, we can identify an ind-scheme 
$Y= ``\ind" Y^\alpha$ with the functor 
$h_Y =\ind \, \Hom(-, Y^\alpha)$ on $\Schb_F$. 
Let $\Affb_F\subset\Schb_F$ 
be the subcategory of affine schemes, dual to $\Algb_F$, the category 
of commutative $F$-algebras. It is well known that $Y$ is determined by the 
restriction 
of $h_Y$ to $\Affb_F$ which can be regarded as a covariant functor on
 $\Algb_F$.

The subcategory of ind-affine ind-schemes over $F$ 
is dual to that of pro-algebras over $F$, i.e., pro-objects $``\pro"A_\a$ 
over 
$\Algb_F$ 
with the structure morphisms beings surjective. Such a pro-algebra can be 
identified with the 
topological algebra $\pro \, A_\a$ with the projective limit topology 
(comp. [Ha]) but we prefer the 
pro-object point of view. For a pro-algebra $A=``\pro" A_\a$ we denote the 
Laurent series algebra  by 
$$A((t)) = \pro  \, A_\a((t)).$$ 
Note that elements of $A((t))$ can be represented by formal series
$\sum_{i=-\infty}^\infty a_i t^i$ with $a_i$ being elements
of the topological algebra $\pro\,  A_\a$. These series are allowed
to be infinite in both directions with the condition that
the $a_i\to 0$ as $i\to -\infty$. 

A morphism $f: Y\to Z$ of ind-schemes will be called an open embedding 
if it is representable by an open embedding of schemes, i.e. 
for any morphism $\phi: S\to  Z $ with $S$ a scheme, the morphism 
$S\times_Z Y \to S$ is a Zariski open embedding of schemes. 
We denote by $Z_{\roman Zar}$ the category of open embeddings 
$Y\to Z$ for a given $Z$. 
 Further, we make $Z_{\roman Zar}$ into a Grothendieck 
site 
 by saying that $\{Y_i\to Y\}$ 
is a covering if for any $\phi, S$ as before the system 
$\{Y_i\times_Z S \to Y\times_Z S\}$ is a Zariski open covering.   
This allows us to speak about $H^i(Z, {\Cal F})$ where $\Cal F$ is a sheaf of 
abelian 
groups on $Z _{\roman Zar}$. 

\proclaim {(2.1) Definition} 
A natural sheaf on the category of $\Schb_F$ is a system 
$\Cal F$ consisting of a Zariski sheaf ${\Cal F}_Y$ given for each scheme $Y$ 
and a morphism of sheaves 
$u_f: f^{-1}({\Cal F}_{Y'})\to {\Cal F}_Y$ given for each 
morphism 
of schemes $f:Y\to Y'$ and satisfying the condition $u_{f\circ g} = u_f\circ 
u_g$ 
for any composable pair of morphisms $f,g$. \endproclaim 

\vskip .2cm 

\noindent {\bf (2.2) Examples.} The sheaves ${\Cal O}_Y$ form a natural 
sheaf. If $\Phi$ is any covariant functor from rings to abelian
groups, then the sheaves $\Phi(\Oc_U)$ (obtained by sheafifying
the presheaves $U\mapsto\Phi(\Oc(U))$ ) form a natural sheaf.
In particular, the sheaves $\Omega^1_Y$ (Kaehler differentials)
 form natural sheaves.

\proclaim{(2.3) Proposition} Let $Z$ be an ind-scheme. 

(a) Any natural sheaf $\Cal F$ gives rise to a sheaf of pro-abelian groups on 
$Z_{\roman Zar}$ 
denoted ${\Cal F}_Z$. In particular, $\Oc_Z$ is a sheaf of pro-algebras on 
$Z_{\roman Zar}$. 

(b) Let $Z$ be represented as $``\ind" \, Z^\alpha$. 
Then we have a natural map 
 $$H^i(Z, {\Cal F}_X) \to \pro \, H^i(Z^\alpha, \Fc_{Z^\alpha})$$ 
which is an isomorphism for $i=0$. 
\endproclaim 

\noindent {\sl Proof:} The Grothendieck site $Z^\alpha_{\roman Zar}$ 
is the category of open sets in $Z^\alpha$ with morphisms being 
embeddings. So whenever $\alpha\leq \alpha'$ we have a functor 
 $Z^{\alpha'}_{\roman Zar}\to Z^\alpha_{\roman Zar}$ 
given by intersecting open sets with $Z^\alpha$. Then, as a category, we have 
$Z_{\roman Zar} = \pro \, Z^\alpha_{\roman Zar}$. Accordingly, the 
Cech complex calculating $H^\bullet(Z, {\Cal F}_X)$ is the 
projective limit of Cech complexes calculating $H^\bullet(Z^\alpha, 
{\Cal F}_{Z^\alpha})$. \qed 

\vskip .2cm 

Note that a morphism $f: Y\to Z$ of ind-schemes gives a moprhism of locally 
ringed 
sites $(Y_{\roman Zar}, \Oc_Y) \to (Z_{\roman Zar}, \Oc_Z)$ and therefore we 
have 
an inverse image map on sheaf cohomology as above. 

\vskip .2cm 

\proclaim {(2.4) Definition}   
Let a natural sheaf $\Cal F$ on the category of 
schemes be given and 
 $Y = ``\ind " Y^\alpha$  be an ind-scheme. An $\Cal F$-gerbe on $Y$ is 
 a system of ${\Cal F}^\alpha$-gerbes ${\Cal 
G}^\alpha$ 
on $Y^\alpha$ together with equivalences of gerbes 
$$\phi_{\alpha\alpha'}: {\Cal G}^\alpha \to i_{\alpha\alpha'}^* ({\Cal 
G}^{\alpha'})$$ 
and isomorphisms of functors 
$$\psi_{\alpha_1\alpha_2\alpha_3}: \phi_{\alpha_1\alpha_3} \Rightarrow 
\phi_{\alpha_1\alpha_2} \circ \phi_{\alpha_2\alpha_3}, \quad 
\alpha_1\leq\alpha_2\leq\alpha_3,$$ 
satisfying the coherence conditions for any 
$\alpha_1\leq\alpha_2\leq\alpha_3\leq\alpha_4$. \endproclaim 

Thus an $\Cal F$-gerbe on $Y$ has a class in $H^2(Y, {\Cal F})$ 
similar to the case of usual schemes. 

\vskip 2cm 

\centerline{\bf 3. The determinantal gerbe of a locally free 
$\Oc_S((t))$-module.} 

\vskip .3cm 

Let $S$ be a scheme and $\Cal E$ be a locally free 
${\Cal O}_S((t))$-module of rank $N$. 
We call a lattice in $\Ec$ a sheaf of $\Oc_S[[t]]$-submodules $\Fc\subset\Ec$ 
such that, Zariski locally on $S$, 
the pair $(\Fc,\Ec)$ is isomorphic to $(\Oc_S[[t]]^N, \Oc_S((t))^N)$. 

\proclaim {(3.1) Proposition} Assume that $S$ is quasicompact. 
Let $\Fc_1, \Fc_2$ are two lattices. 

(a) 
Zariski locally on $S$, there exist $a,b\in\ZZ$ such that 
$$t^a \Fc_1 \subset \Fc_2\subset t^b \Fc_1.$$ 

(b) 
If $\Fc_1\subset\Fc_2$ then $\Fc_2/\Fc_1$ is locally 
free over $\Oc_S$ of finite rank. 
\endproclaim 

\noindent {\sl Proof:} 
(a) By quasicompactness of $S$ it is 
enough to prove the statement locally, in particular, to 
assume that $\Ec=\Oc_S((t))^N$, that $\Fc_2=\Oc_S[[t]]^N$ 
and that $\Fc_1 = A\cdot \Fc_2$ where $A\in\GL_N(\Oc(S)((t)))$. 
Our statement follows by taking  $a$ to be the negative 
of the maximal order of poles of the coefficients of $A$ 
and $b$ to be the maximal order of poles of the coefficients of 
$A^{-1}$. 

(b) 
Again, it is enough to work locally and assume that $S=\Spec(R)$ is affine. 
Then the existence of $A$ implies that for any affine $U\subset S$ 
the $\Oc(U)[[t]]$-submodule 
$\Fc_i(U)\subset\Ec(U)$ is coprojective 
over $\Oc(U)$, i.e. $\Ec(U)/\Fc_i(U)$ is projective. 
It follows that the sequence of $\Oc_S$-modules 
$$0\to\Fc_2/\Fc_1\to\Ec/\Fc_1\to\Ec/\Fc_2\to 0$$ 
splits over $\Oc_S$, so $\Fc_2/\Fc_1$ is locally a direct summand 
of a projective $\Oc_S$-module, hence it is itself projective. 
The finiteness of the rank follows from (a). 
\qed 

\subhead{(3.2) The twisted affine Grassmannian} 
\endsubhead 

Consider the contravariant functor  $\gamma: \Schb_S\to\Setsb$ 
taking $p\,:\,T\to S$ to the set of 
$\Oc_T[[t]]$-lattices $\Fc\subset p^*\Ec$. Here 
$$p^*\Ec = p^{-1}\Ec \otimes_{p^{-1}\Oc_S((t))} \Oc_T((t))$$ 
is the direct image in the sense of morphisms of ringed spaces 
$$(S, \Oc_S((t)))\to (T,\Oc_T((t))).$$

\proclaim {(3.2.1) Proposition} The functor $\gamma$ 
is represented by an ind-proper formally smooth ind-scheme 
over  $S$, denoted by $\Grc(S, \Ec)$. 
\endproclaim 

\noindent {\sl Proof:} It is enough to prove the statement Zariski 
locally on $S$. In particular, we can assume that $S = \Spec (R)$ 
is affine and that 
 $\Ec =\Oc_S((t))^N$ is free. For any $a<b \in \ZZ$ let 
$\gamma_{a,b}$ be the sub-functor of $\gamma$ taking $T$ to 
the set of lattices $$\Fc\subset p^*\Ec = \Oc_T((t))^N$$ 
that are contained between 
$t^a \Oc_T[[t]]^N$ and $t^b\Oc_T[[t]]^b$. Let also $\xi_{a,b}$ be the 
functor taking $T$ to the set of $\Oc_S$-submodules contained between 
$t^a \Oc_T[[t]]^N$ and $t^b\Oc_T[[t]]^b$. Then, clearly, 
$\xi_{a,b}$ is represented by the relative Grassmannian 
of subspaces (of all dimensions) in the vector bundle 
$t^b\Oc_S[[t]]^N/t^a\Oc_S[[t]]^N$ and $\gamma_{a,b}$ 
is defined inside this relative Grassmannian by closed conditions 
of being an $\Oc_S[[t]]$-submodule. So $\gamma_{a,b}$ 
is represented by a proper scheme $\Grc_{a,b}(S,\Ec)$ 
of finite type over $S$. Finally, as 
$\gamma = \ind _{a,b} \gamma_{a,b}$ in the category of functors, 
we conclude that it is representable by the ind-scheme 
$\Grc(S, \Ec) = ``\ind ''
)_{a,b} \Grc_{a,b}(S, \Ec)$. \qed

\subhead{(3.3) The ${\Cal O}^\times$-groupoid structure}\endsubhead 

Proposition 3.1 can be reformulated as follows.

\proclaim{(3.3.1) Lemma} Given an open set $U\subset S$ and two $U$-points 
$\Fc_1, \Fc_2$ of $\Grc(S, \Ec)$, we can find Zariski locally on $U$ 
a third point $\Fc'$ such that as a submodule $\Fc'$ is contained 
in each of $\Fc_i$ 
and the quotients $\Fc_i/\Fc'$ are 
free over $\Oc_U$ of finite rank. 
\endproclaim 

We now define 
$$(\Fc_1|\Fc_2) = \Wedge^{\max}(\Fc_1/\Fc')\otimes\Wedge^{\max} 
(\Fc_2/\Fc')^{-1}.$$ 
This is a line bundle that is independent (up to a unique isomorphism) 
on the choice of $\Fc'$ and thus it is well defined by gluing for any 
two $U$-points $\Fc_1, \Fc_2$. It is also clear that for three $U$-points we 
have a canonical isomorphism 
$$(\Fc_1|\Fc_2)\otimes (\Fc_2|\Fc_3) \to (\Fc_1|\Fc_3).$$ 

\proclaim{(3.3.2) Lemma} The $\Oc_U^\times$-torsors corresponding to the 
line bundles 
$(\Fc_1|\Fc_2)$ make $\Grc(S, \Ec)$ into the sheaf of objects of a sheaf of 
$\Oc^\times$-groupoids on $S$. 
This sheaf of groupoids is locally connected and hence 
gives rise to an $\Oc_S^\times$-gerbe. 
\endproclaim 

We denote this gerbe by $\Det(\Ec)$. 
Note that $[\Det(\Ec)]\in H^2(S,\Oc_S^\times)$. 

\vskip 2cm 

\centerline {\bf 4. Loop ind-schemes and the evaluation maps.} 

\vskip .3cm 

\subhead{(4.1) The restriction of scalars} 
\endsubhead 
From now all all rings will be assumed to contain the field $\CC$ of complex
numbers and by (ind-)schemes we will mean (ind-)schemes over $\CC$. 

We define a functor 
$$\Rc=\Rc_\CC^{\CC((t))}: \Affb_{\CC((t))}\to\Ischb,$$ 
called the restriction of scalars (from $\CC((t))$ to $\CC$). 
Given an affine scheme $Y$ over $\CC((t))$, the ind-scheme 
$\Rc Y$ represents the following functor on $\Algb$ : 
$$\Hom_{\Ischb}(\Spec A,\Rc Y)= 
\Hom_{\Schb_{\CC((t))}}(\Spec A((t)),\Rc Y).$$ 
It is well-known that this functor is indeed represented by an ind-affine 
ind-scheme. 
Compare \cite{D, Sect.~6.3.4}, where $\Rc Y$ is denoted by $\Lc Y$. 

\proclaim{\bf (4.1.2) Proposition} 
(a) The functor $\Rc$ commutes with finite projective limits, in particular, 
with fiber products. 

(b) The functor $\Rc$ also commutes with disjoint unions. 
\endproclaim 

\noindent{\sl Proof :} 
(a) This follows directly from (4.1.1) as, for any $S$, 
$$\Hom(S,\pro \, Z_i)=\pro\, \Hom(S,Z_i).$$ 

(b) Obvious. 
\qed

\subhead{(4.2) The evaluation map} \endsubhead 
  
Let $Y\in\Affb_{\CC((t))}$. 
We define a morphism of ringed Grothendieck sites 
$$\epsilon:\biggl((\Rc Y)_\Zar,\Oc_{\Rc Y}((t))\biggr)\to(Y_\Zar,\Oc_Y)$$ 
called the evaluation map. The underlying morphism of sites 
$$\eps_\sharp:(\Rc Y)_\Zar\to Y_\Zar$$ 
is defined by putting, for a Zariski open $U\subset Y$ 
$$\epsilon^{-1}(U)=\Rc U\to\Rc Y.$$ 

\proclaim{(4.2.1) Proposition} 
(a) $\Rc U$ is Zariski open in $\Rc Y$. 

(b) The correspondence $U\mapsto \Rc U$ 
commutes with disjoint unions and fiber products. 

(c) If $U=\bigcup_\alpha U_\alpha$ is a covering in $Y_\Zar$, then 
$\{\Rc U_\alpha\}$ form a covering of $\Rc U$ in 
$(\Rc Y)_\Zar$. 
\endproclaim 

\noindent{\sl Proof :} 
(a) Let $Y=\Spec R$, $R\in\Algb_{\CC((t))}$. 
A basis of open sets in $Y_\Zar$ is formed by the "principal affine sets" 
$$U_f=\{x\in Y: f(x)\neq 0\},\quad f\in R.$$ 
It is enough to prove that $\Rc U_f$ is open in $\Rc Y$. 
Let $y$ be a $\CC$-point of $\Rc Y$ and $y(t)$ be the corresponding 
$\CC((t))$-point of $Y$. 
Thern $y\in\Rc U_f$ if and only if $y(t)\in U_f$, i.e. 
$f(y(t))\neq 0$ as an element of $\CC((t))$. 
Writting $f(y(t))=\sum_{n\in\ZZ}a_n t^n$, 
we realize each $a_n$ as a global section 
of $\Oc_{\Rc Y}$. 
So 
$$U_{a_n}=\{y\in\Rc Y: a_n(y)\neq 0\}$$ 
is Zariski open in $\Rc Y$ and 
$$\Rc U_f=\bigcup_{n\in\ZZ}U_{a_n}$$ 
is open as well. 

(b) Follows from Proposition (4.1.2). 

(c) If $U=\bigcup_\alpha U_\alpha$, then 
for any field $F$ we have $U(F)=\bigcup_\alpha U_\alpha(F)$. 
Taking $F=\CC((t))$, we get 
$$(\Rc U)(\CC)=U\biggl(\CC((t))\biggr)=
\bigcup_\alpha U_\alpha \bigl(\CC((t))\bigr)= 
\bigcup_\alpha(\Rc U_\alpha)(\CC).$$ 
\qed 

\vskip3mm 

Let $S=\Spec A$ be an affine scheme over $\CC$, and 
$f:S\to\Rc Y$ be a morphism. 
Let $\tilde f:\Spec A((t))\to Y$ 
be the morphism of $\CC((t))$-schemes 
corresponding to $f$ by (4.1.1). 
Then $\Fc$ induces a ring homomorphism 
$$\tilde f^*:\Gamma(Y,\Oc_Y)\to A((t))=\Gamma\bigl(S, \Oc_S((t))\bigr).$$ 
By performing the same construction for affine open sets 
$U\subset Y$ (forming a basis of topology of $Y_\Zar$) we get 
a morphism of sheaves of rings on $S$ 
$$\epsilon^\flat_{S,f}:f^{-1}\epsilon_\sharp^{-1}\Oc_Y\to\Oc_S((t)).$$ 
These morphisms, taken for all $S$, $f$ as above, give rise to 
a morphism of sheaves of rings on $\Rc Y$ 
$$\epsilon^\flat:\epsilon_\sharp^{-1}\Oc_Y\to\Oc_{\Rc Y}((t)).$$ 
This finishes the definition of $\epsilon$. 

\subhead{(4.3) The transgression map} \endsubhead 
 We now 
discuss an algebraic construction which is an analog of 
$\int_{\Sigma/B}$ in (0.3.3). 

Let $A$ be a commutative ring. We denote by $K_j(A)$ the $j$th 
Quillen K-group, so $K_*(A)$ is a graded commutative ring. 
Let also $K_j^M(A)$ be the $j$th Milnor K-group. Recall [EM] 
that $K_*^M(A)$ is defined as the graded commutative ring 
generated by symbols $l(a), a\in A^\times$ in degree 1 which are subject 
to the relations: 
$$l(ab) = l(a) + l(b), \quad a,b\in A^\times;$$ 
$$l(a) l(1-a) =0, \quad a, 1-a\in A^\times;$$ 
$$l(a) l(-a) =0, \quad a\in A^\times.$$ 
In particular,  $K_0^M(A) = \Zb$ and $K_1^M(A) = A^\times$. 
One denotes the element $l(a_1)... l(a_n)\in K_n^M(A)$ 
by $\{a_1, ..., a_n\}$. There is a canonical ring homomorphism 
$$\lambda_*: K_*^M(A)\to K_*(A),\leqno (4.3.1)$$ 
taking $a\in A^\times$ to its class in $K_1(A)$. The following fact 
follows from the results of [VdK]. 

\proclaim{(4.3.2) Theorem} If $A$ is a local ring with an infinite 
residue field, then $\lambda_i: K_i^M(A)\to K_i(A)$ 
is an isomorphism for $i\leq 2$. 
\endproclaim 

Let $S$ be a scheme over $\Cb$. By sheafifying the presheaf 
 $U\mapsto K_j^M 
(\Oc(U))$ we get a sheaf $K_j^M(\Oc_S)$. Since the 
stalks of $\Oc_S$ are local $\Cb$-algebras, the above theorem 
implies: 

\proclaim{(4.3.3) Corollary} The sheaves $K_j^M(\Oc_S)$ and 
$K_j(\Oc_S)$ are identified for $j\leq 2$. 
\endproclaim 

We now recall the definition of  the homomorphism 
$$\partial: K_2^M A((t))\to A^\times=K_1^M A,\leqno (4.3.4)$$ 
called the Contou-Carrere symbol [CC]. To do this, it suffices 
to define $\partial \{a(t), b(t)\}\in A^\ast$ for any two units $a(t),b(t) 
\in A((t))^\times$ and show that it satisfies the relations for $K_2^M$ 
above. We recall the following lemma, see  [CC] (1.3) and [De2] (2.9)

\proclaim {(4.3.4) Lemma} (a) Assume that $\Spec (A)$ is connected,
i.e., $A$ is not a product of rings. Then
any invertible $a(t)=\sum_{i\gg -\infty}^\infty a_it^i\in A((t))^\times$
satisfies the following property: There exist $n = {\text ord}(a)\in \ZZ$
 such
that $a_n$ is invertible and $a_i, i<n$ are nilpotent. 

(b) The correspondence $a\mapsto {\text ord}(a)$ is a homomorphism
$A((t))^\times \to \ZZ$.

(c) If $A$ is arbitrary and $a(t)\in A((t))^\times$, then there
is a decomposition of $A$ into a finite product of rings such that
the description from (a) applies to each factor.

\endproclaim

This can be reformulated as a statement describing the ind-scheme
$\GL_1((t)) = \Rc\bigl(\GL_1/\CC((t))\bigr)$. By definition,
an $A$-point of $\GL_1((t))$ is an element of $\GL_1\bigl( A((t))\bigr) = A((t))^\times$. 

\proclaim{(4.3.5) Corollary}We have an identification of ind-schemes:
$$\GL_1((t))  =  \ZZ \times \CC^\times\times \Spec \, \CC[a_1, a_2, ...] \times \Spf
\, \CC [[ a_{-1}, a_{-2}, ...]].$$
Here $\Spf
\, \CC [[ a_{-1}, a_{-2}, ...]]$ is, by definition, the ind-scheme
$$``\ind''_\epsilon \Spec \CC[ a_{-1}, a_{-2}, ...]\biggl/ \bigl ( a_i^{\epsilon_i+1} = 0, 
\, i<0\bigr),$$
where $\epsilon$ runs over sequences $\epsilon=(\epsilon_{-1}, \epsilon_{-2}, ...)$,
 $\epsilon_i\in \ZZ_+$, almost
all $\epsilon_i=0$.
\endproclaim

The above identification corresponds to writing an element of $A((t))$ in (4.3.4)(a)
as a product
$$a(t) = t^n \cdot a_0 \cdot \biggl(1+\sum_{i=1}^\infty a_i t^i\biggr)
 \cdot \biggl(1+\sum_{i=-\infty}^{-1} a_i t^i
\biggr),$$
with $a_0\in A^\times$, $a_i, i>0$ arbitrary and $a_i, i<0$ nilpotent. 

\vskip .2cm

Let us now describe the Contou-Carrere symbol.

In virtue of (4.3.4) (c) we will assume in the remainder of this
section that $A$ is not a product. 
Let $a(t), b(t)\in A((t))^*$ be two units with ${\text ord}(a(t)) = n,
{\text ord}(b(t)) = m$. The Contou-Carrere symbol of $a(t)$ and $b(t)$
is defined by (see [De2](2.9)):
$$\partial\{a(t), b(t)\} := (-1)^{mn} b_m^{-n} \exp\left( \Res\left(
{da\over a} \cdot \log\left({b(t)\over b_mt^m}\right)\right)\right).
\leqno (4.3.6)$$
Here $b(t)/b_mt^m$ is the sum of 1 and a topologically nilpotent
element so its logarithm is a well defined element of $A((t))$.

\vskip .3cm

\noindent{\bf (4.3.7) Remark.} For any commutative ring $A$
one can obtain the map $\partial^Q_i: K_i A((t))\to K_{i-1} A$
on Quillen K-functors using the localization theorem for
singular varieties ([Sr], Ch. 9). It seems that the image of
$\partial^Q_2$ on elements $\lambda_2\{a(t), b(t)\}\in K_2 A((t))$
for $a(t), b(t)\in A((t))^\times$ is in fact equal to the
Contou-Carrere symbol although we could not find a reference for
this in the literature. 

\vskip .3cm

Let now $Y$ be a scheme over $\CC ((t))$. 
For any sheaf $\Ec$ of $\Oc_Y$-modules we have a sheaf 
$$\Rc\Ec=\epsilon^*\Ec= 
\epsilon_\sharp^{-1}\Ec\otimes_{\epsilon_\sharp^{-1}\Oc_Y}\Oc_{\Rc Y}((t))$$ 
of $\Oc_{\Rc Y}((t))$-modules. In particular, 
$\Rc \Oc_Y=\Oc_{\Rc Y}((t))$.

The morphism of ringed spaces $\epsilon$ induces homomorphisms 
of Abelian groups 
$$\epsilon^\circ:H^i(Y,K^M_j(\Oc_Y))\to H^i\biggl(\Rc Y,K^M_j\bigl(\Oc_{\Rc Y}((t))\bigl)\biggl).$$

Sheafifying the Contou-Carrere symbol  on $\Rc Y$, we get a morphism
 of sheaves of Abelian groups 
$$\partial: K_2^M(\Oc_{\Rc Y}((t)))\to K_{1}^M(\Oc_{\Rc Y})=
 \Oc_{\Rc Y}^\times.$$ 
Composing it with $\epsilon^\circ$ we get, for any $i\ge 0$, 
a homomorphism 
$$\tau=\tau_{i}: H^i(Y,K_2(\Oc_Y))\to H^i(\Rc Y,\Oc_{\Rc Y}^\times)
\leqno (4.3.8)$$ 
which we call the transgression.

\vskip 2cm 

\centerline {\bf 5. Chern classes and the local Riemann-Roch theorem.} 

\vskip .3cm 

\subhead{(5.1) Reminder on simplicial geometry} \endsubhead 

Recall that a simplicial object in a category $\Cb$ is  a familly 
$X_\bullet=(X_n)_{n\ge 0}$ of objects of $\Cb$ together with 
face and degeneracy morphisms 
$$d_i:X_n\to X_{n-1},\quad 
s_i:X_n\to X_{n+1},\quad 
i=0,1,...n$$ 
satisfying the standard identities. 

\vskip3mm 

\noindent{\bf (5.1.1) Examples.} 
(a) 
Let $F$ be a field and $\Cb=\Schb_F$. 
Let $S$ be a scheme and $U=\{U_\alpha\}_\alpha$ be a 
Zariski open covering. The nerve $\Nc_\bullet U$ is a simplicial scheme 
with 
$$\Nc_nU=\coprod_{\a_0,\a_1,..a_n\in A}(U_{\alpha_0}\cap\cdots
 U_{\alpha_n}).$$ 

(b) 
Let $\Cb=\Ischb_F$ and $G$ be a group ind-scheme over $F$. 
The classifying space of $G$ is the simplicial ind-scheme 
$B_\bullet G$ with $B_nG=G^n$. 

\vskip3mm 

Let $X_\bullet$ be a simplicial topological space. 
A sheaf of Abelian groups $\Fc_\bullet$ on $X_\bullet$ 
is a family of sheaves $\Fc_n$ on $X_n$ together with morphisms 
$$d_i^{-1}\Fc_{n-1}\to\Fc_n,\quad 
s_i^{-1}\Fc_{n+1}\to\Fc_n,$$ 
compatible with the identities among $d_i$, $s_i$. 
The groups 
$$\HH^i(X_\bullet,\Fc_\bullet)$$ 
are defined as the cohomology groups of the double complex formed by the Cech 
complexes of $\Fc_n$ on $X_n$ for $n\ge 0$, see \cite{G}, 
so that one has a spectral sequence 
$$E_2=H^i(X_j,\Fc_j)\Rightarrow\HH^{i-j}(X_\bullet,\Fc_\bullet).$$ 

\vskip1mm 

\noindent{\bf (5.1.2) Examples.} 
(a) 
If $X_\bullet$ is a simplicial ind-scheme, then 
we have the simplicial topological space $(X_\bullet)_\Zar$, and 
$\Oc=\Oc_{X_\bullet}=(\Oc_{X_n})$ is a sheaf of rings 
on $(X_\bullet)_\Zar$. 
Therefore 
$$K_j(\Oc)=(K_j(\Oc_{X_n}))$$ 
is a sheaf of Abelian groups on 
$(X_\bullet)_\Zar$. 
Similarly we have the sheaf 
on $(X_\bullet)_\Zar$ 
$$K_j(\Oc((t)))=\biggl(K_j\bigl(\Oc_{X_n}((t))\bigl)\biggl).$$ 

(b) 
If $S$ is a scheme, $U=\{U_\alpha\}_\alpha$ is a Zariski 
open cover, and $\Fc$ is a sheaf on $S_\Zar$, 
we get a simplicial sheaf $\Fc_\bullet$ on $\Nc_\bullet U$ 
with $$\Fc_n|_{U_{\alpha_0}\cap\cdots U_{\alpha_n}}= 
\Fc|_{U_{\alpha_0}\cap\cdots U_{\alpha_n}}.$$ 
The Mayer-Vietoris property for sheaf cohomology implies that 
$$\HH^i(\Nc_\bullet U,\Fc_\bullet)=H^i(S,\Fc).$$ 

\subhead{(5.2) The Chern classes of Goncharov} \endsubhead 

Let $F$ be a field and $\GL_N/F$ be the algebraic group $\GL_N$ 
considered as a scheme over $F$. 
Then $B_\bullet(\GL_N/F)$ is a simplicial $F$-scheme. 
Goncharov [Go]  has constructed classes 
$$c_i\in\HH^i(B_\bullet(\GL_N/F), K^M_i(\Oc))$$ 
which we will call  the universal Chern classes. 
They lift the Chern classes of Gillet [Gi]
to Milnor K-theory. In this paper we will be working only with
$c_i$ for $i=1,2$ when, by Corollary 4.3.3, there is no
difference between $K_i$ and $K_i^M$ for sheaves of local 
$\Cb$-algebras.

\vskip3mm 

\noindent{\bf (5.2.1) Example.} 
Let $S$ be a smooth scheme over $F$ and $\Ec$ 
a locally free sheaf of $\Oc_S$-modules of rank $N$. 
Choosing a Zariski open covering $U=\{U_\alpha\}_\alpha$ 
of $S$ such that $\Ec$ is free over each $U_\alpha$, 
and a system of trivializations 
$$\varphi=(\varphi_\alpha:\Ec|_{U_\alpha}\to\Oc_{U_\alpha}^N)$$ 
one gets the system of transition functions 
$$\varphi_{\alpha,\beta}:U_\alpha\cap U_\beta\to\GL_N.$$ 
These functions give a morphism of simplicial schemes 
$$\tilde\varphi:\Nc_\bullet U\to B_\bullet(\GL_N/F).$$ 
The inverse image 
$$\widetilde{\varphi}^*c_i\in\HH^i(\Nc_\bullet U,K_i^M(\Oc))
=H^i(S,K_i^M(\Oc_S))$$ 
is nothing but $c_i(E)$, the $K$-theoretic Chern class of $\Ec$. 
To be precise, for a Noetherian regular scheme $S$ we have
$H^i\bigl(S, K_i(\Oc_S)\bigr)= CH^i(S)$ is the Chow group
of cycles of codimension $i$, see [Sr]. The image of
$\widetilde{\varphi}^* c_i$ under the map
$\lambda_i$ of (4.3.1) is the $i$th Chern class with values in
the Chow group.

\vskip3mm 

Let us now take $F=\CC((t))$. 
The group scheme $\GL_N/\CC((t))$ gives the group ind-scheme 
$$\GL_N((t)):=\Rc\bigl(\GL_N/\CC((t))\bigl)$$ 
over $\CC$. The evaluation map from (4.3) 
gives a morphism of ringed simplicial topological spaces 
$$\epsilon: 
(B_\bullet\GL_N((t)),\Oc((t)))\to 
\biggl(B_\bullet\bigl(\GL_N/\CC((t))\bigl),\Oc\biggl).\leqno(5.2.2)$$ 
Therefore we have the classes 
$$c_i^{((t))}=\epsilon^\circ c_i\in 
\HH^i\biggl(B_\bullet\bigl(\GL_N((t))), K_i^M(\Oc((t))\bigl)\biggl).$$ 
Let now $S$ be a scheme over $\CC$ and $\Ec$ 
be a locally free sheaf of $\Oc_S((t))$-modules of rank $N$. 
Trivializing $\Ec$ over open sets from a covering $U$, as in Example (5.2.1), 
we get a morphism of simplicial ind-schemes over $\CC$ 
$$\tilde\varphi:\Nc_\bullet U\to B_\bullet\GL_N((t))\leqno(5.2.3)$$ 
which induces the classes 
$$c_i^{((t))}(\Ec):=\tilde\varphi^*c_i^{((t))}\in 
H^i\biggl(S,K_i^M\bigl(\Oc_S((t))\bigl)\biggl).$$ 

\subhead{(5.3) The local Riemann-Roch Theorem} \endsubhead 

Let $S$ be a scheme, $\Ec$ a locally free sheaf of $\Oc_S((t))$-modules 
of rank $N$. 
Then we have the second Chern character class 
$$ ch_2^{((t))}(\Ec)={1\over 2} \bigl(c_1^{((t))}(\Ec)\bigr)^2-
c_2^{((t))}(\Ec) 
\in H^2\biggl(S,K_2^M\Oc_S((t))\biggr)\otimes \ZZ \bigl[{1\over 2}\bigr].$$ 
On the other hand we have the determinantal class 
$$[\Det(\Ec)]\in H^2(S,O^\times_S).$$ 
The Contou-Carrere  symbol 
gives a homomorphism of sheaves 
$$\partial:K_2\bigl(\Oc_S((t))\bigr)\to\Oc_S^\times.$$ 
It yields a map 
$$\partial:H^2\biggl( S,K_2\bigl( \Oc_S((t))\bigr) \biggr) \to H^2(S,\Oc_S^\times).$$ 

\proclaim{(5.3.1) Theorem} 
We have $$[\Det(\Ec)]=\partial( ch_2^{((t))}(\Ec)) \quad {\text in}\quad 
H^2(S, \Oc_S^\times)\otimes \ZZ \bigl[{1\over 2}\bigr].$$ 
\endproclaim 

To prove the theorem, we realize $[\Det(\Ec)]$ as 
$$[\Det(\Ec)]=\tilde\varphi^*[\Det]\leqno(5.3.2)$$ 
for a universal class 
$$[\Det]\in\HH^2(B_\bullet\GL_N((t)),\Oc^\times).\leqno(5.3.3)$$ 
Here $\tilde\varphi$ is as in (5.2.3). 
Then the theorem reduces to the statement about the universal classes : 
$$[\Det]=\partial({1\over 2}(c_1^{((t))})^2-c_2^{((t))})\in 
\HH^2(B_\bullet\GL_N((t)),\Oc^\times)\otimes \ZZ\bigl[{1\over 2}\bigr].$$ 
We start with a standard lemma. 

\proclaim{(5.3.5) Lemma} 
Let $G$ be a group ind-scheme over $\CC$. 
Then $\HH^2(B_\bullet G,\Oc^\times)$ is identified 
with the group of central extensions 
$$1\to\CC^\times\to\widetilde {G}\to G\to 1$$ 
(in the category of group ind-schemes). 
\endproclaim 

Recall now the definition of the determinantal central extension 
$$1\to\CC^\times\to\widetilde{\GL}_N((t))\to\GL_N((t))\to 1.$$ 
The group ind-scheme $\widetilde{\GL_N}((t))$ represents the following functor 
on $\Algb$ : 
$$\Hom(\Spec(A),\widetilde {\GL}_N((t)))=$$ 
$$\biggl\{(g,u): g\in\GL_NA((t)),\quad u\in\det(gA[[t]]^N:A[[t]]^N)\
 \roman{invertible}\biggr\}. 
\leqno(5.3.6)$$ 
Here ``invertible" means``an element of a projective $A$-module of rank one 
which does not vanish over any prime ideal". 
We now define $[\Det]$ in (5.3.3) to be the class of 
$\widetilde{\GL}_N((t))$. 
It is clear from the definition of 
$[\Det(\Ec)]$ via the determinantal gerbe that the equality (5.3.2) holds. 
We concentrate therefore on the proof of (5.3.4), which, by (5.3.5), is a 
statement comparing two central extensions of $\GL_N((t))$. 

\subhead{(5.4) Proof of Theorem 5.3.1} \endsubhead 

Let $\widetilde {G}$ be the central extension of $\GL_N((t))$ 
corresponding to the class $2[\Det]-\partial(2\ch_2^{((t))})$. 
We will prove that $\widetilde {G}$ is trivial. 

\proclaim{(5.4.1) Proposition} 
We have that $\widetilde {G}$ is trivial on the sub-ind-group scheme 
$T((t))$ where $T\subset\GL_N$ is the subgroup of diagonal matrices. 
\endproclaim 

\noindent{\sl Proof :} 
To emphasize the dependence of our classes on $N$, let us denote them 
by $[\Det]_N$ and $2\ch_{2,N}^{((t))}$. 
Notice that both these classes are additive with respect to the block diagonal 
embeddings 
$$\varphi_{N_1,N_2}:\GL_{N_1}((t))\times\GL_{N_2}((t))\to\GL_{N_1+N_2}((t)),$$ 
i.e. 
$$\varphi_{N_1,N_2}^*([\Det_{N_1+N_2}])= 
[\Det_{N_1}]+[\Det_{N_2}]$$ 
and similarly for $\ch_2^{((t))}$. 
Indeed, for $[\Det]$ it follows right away from 
(5.3.6) and for $2\ch_2^{((t))}$ by the standard 
additivity property of the Chern character. 
Therefore it is enough to establish the lemma for $N=1$. 
Then $2\ch_2^{((t))}=(c_1^{((t))})^2$. 
Let us first identify $c_1^{((t))}$. 
Notice that the hypercohomology spectral sequence 
gives a homomorphism 
$$\xi:\Hom\bigl(\GL_1((t)),\GL_1((t))\bigl)\to
\HH^1(B_\bullet\GL_1((t)),\Oc((t))^\times)$$ where $\Hom$ 
on the left is the set of homomorphisms of group ind-schemes. 
The next lemma is clear.

\proclaim{(5.4.2) Lemma} 
We have that $c_1^{((t))}=\xi(\Id)$. 
\endproclaim\qed 

Lemma 5.4.2 implies that $\partial(c_1^{((t))})^2$ is the image of the 
tame symbol map 
$$\GL_1((t))\times\GL_1((t))\to\CC^\times$$ 
under the homomorphism 
$$Z^2(\GL_1((t)),\CC^\times)\to\HH^2(B_\bullet\GL_1((t)),\Oc^\times)$$ 
where $Z^2$ is the set of 2-cocycles that are morphisms of ind-schemes. 
To finish the proof of (5.4.1), it suffices to establish  the following. 

\proclaim {(5.4.3) Lemma} In the group 
$\HH^2(B_\bullet \GL_1((t)), \Oc^*)\otimes \ZZ\bigl[{1\over 2}]$ we have 
$$2[\Dc et] = \tau = \partial((c_1^{((t))})^2).$$ 
\endproclaim 

\noindent{\sl Proof:} We start with a general remark. Let $E, A$ be abelian
groups and $\widetilde{E}$ be a central extension
$$1\to A\to \widetilde{E}\to E\to 1.$$
As any central extension of any  group (abelian or not),
 it corresponds to a class
in $H^2(E, A)$. The well known procedure to represent  this class by a cocycle
is as follows. One takes
for each $x\in E$ a lifting $\widetilde{x}\in\widetilde{E}$ and
 writes the cocycle
$$\gamma(x,y) = \widetilde{x} \widetilde{y}.(\widetilde{xy})^{-1}.$$
On the other hand,  since $E$ is abelian, one can form the
commutator pairing
$$c(x,y) = \widetilde{x} \widetilde{y} \widetilde{x}^{-1}\widetilde{y}^{-1},$$
which, as any bilinear pairing on any abelian group,
 can also be considered as a 2-cocycle. 
It is clear that 
$$c(x,y) = \gamma(x,y) \gamma(y,x)^{-1}.$$
On the other hand, $\gamma(x,y)\gamma(y,x)$ is the coboundary
of the 1-cochain $x\mapsto\gamma(x,x)$. This means that
the class if $c(x,y)$ is equal to twice the class of $\gamma(x,y)$. 

\vskip .1cm

We now apply the above to the case where $E$ is the group scheme
$GL_1((t))$, $A= \CC^*$ and $\widetilde{E}$ is the central extension
corresponding to the determinantal gerbe. To be precise, we take
a commutative algebra $A$ and apply the above to  groups of $A$-points. 
It follows from the results of [AP] that for $a(t), b(t)\in A((t))^*$
 the commutator pairing $c(a(t), b(t))$ is equal to 
$(-1)^{{\text ord}(a) {\text ord}(b)}$ times the Contou-Carrere
symbol of $a(t), b(t)$. To be precise, in [AP] the authors treated the
case of an artinian local ring $A$. However, both the commutator pairing and
the Contou-Carrere symbol are regular functions, i.e., morphisms of ind-schemes
$$\GL_1((t))\times\GL_1((t))\to \CC.$$
It follows from the explicit description of $\GL_1((t))$ as an ind-scheme
(Corollary 4.3.4) that two such morphisms are equal if and only if they
are equal on the set of points with values in any artinian local ring. 

It follows that after inverting 2,
the cohomology class describing the extension is one half of the
class of the Contou-Carrere symbol, i.e.,
$2[\Dc et] = \partial((c_1^{((t))})^2)$ as claimed in Lemma 5.4.3. 
This also finishes the proof of
Proposition 5.4.1.

\proclaim{(5.4.4) Lemma} 
The extension $\widetilde {G}$ is trivial on the subgroup scheme $\GL_N[[t]]$. 
\endproclaim 

\noindent{\sl Proof :} 
For $[\Det]$ the triviality on $\GL_N[[t]]$ follows from (5.3.6) and for 
$\partial(2\ch_2^{((t))})$ from the fact that $\partial$ is trivial 
on $K_2(A[[t]])$ for any ring $A$. 
\qed 

\proclaim{(5.4.5) Lemma} 
The extension $\widetilde {G}$ is trivial on the subgroup ind-scheme $\SL_N((t))$. 
\endproclaim 

\noindent{\sl Proof :} 
Let $\Gc=\SL_N((t))/\SL_N[[t]]$ be the affine Grassmannian for $\SL_N$. 
Because of (5.4.4) the extension $\tilde G|_{\SL_N((t))}$ 
(viewed as a multiplicative $\Oc^\times$-torsor on 
$\SL_N((t))$) descends to an $\Oc^\times$-torsor on $\Gc$, i.e. 
to a line bundle which we denote by $\Lc$. 
Recall \cite{Ku, Proposition 13.2.19} that $\Pic(\Gc)=\ZZ$. 
Let $T_{\SL}\subset T$ be the intersection $T\cap\SL_N$. 
By Lemma 5.4.1 the line bundle $\Lc$ is equivariant with respect to 
$T_\SL((t))$. However, it is known that any nontrivial line bundle 
on $\Gc$ is equivariant with respect to a nontrivial central extension 
of $\SL_N((t))$, and this extension remains nontrivial on 
$T_\SL((t))$. So $T_\SL((t))$-equivariance implies that $\Lc$ 
is trivial. This implies that the projection 
$$\widetilde {G}|_{\SL_N((t))}\to\SL_N((t))$$ 
splits as a morphism of ind-schemes and so 
$\widetilde {G}|_{\SL_N((t))}$ is given by a 2-cocycle 
$$\eta:\SL_N((t))\times\SL_N((t))\to\CC^\times$$ 
which is a morphism of ind-schemes. 
The fact that $\eta$ is trivial follows from the next lemma. 

\proclaim{(5.4.6) Lemma} 
We have that $\Gamma(\SL_N((t)),\Oc^\times)=\CC^\times$. 
\endproclaim 

\noindent{\sl Proof :} 
As well-known, any unimodular matrix over any field, in 
particular, over the field $\CC((t))$ 
can be factored as a product of elementary matrices 
$e_{ij}(a)$, $a\in\CC((t))$, $i,j=1,2,...N$, $i\neq j$. 
This means that there is a sequence of pairs of indices $(i_\nu,j_\nu)$, 
$\nu=1,2,...M$ such that every matrix as above can be written as 
$$e_{i_1,j_1}(a_1)e_{i_2,j_2}(a_2)...e_{i_M,j_M}(a_M),\quad 
a_\nu\in\CC((t)).$$ 
Let now $\underline{\CC((t))}=\Rc(\AA^1/{\CC((t))})$ be the 
ind-scheme whose set of $\CC$-points is $\CC((t))$. 
Asit is an inductive limit of affine spaces over $\CC$, we have 
$$\Gamma(\underline{\CC((t))},\Oc^\times)=\CC^\times.$$ 
On the other hand, the construction above gives a surjective 
morphism of ind-schemes 
$$\underline{\CC((t))}^M\to\SL_N((t)).$$ 
A nonvanishing function $f$ on $\SL_N((t))$ gives then a nonvanishing 
function $\tilde f$ on $\underline{\CC((t))}^M$ which must be constant. 
So $f$ is constant as well. 
Lemmas 5.4.6 and 5.4.5 are proved. 
\qed 

\vskip3mm 

Now we prove the triviality of $\widetilde {G}$ on the whole $\GL_N((t))$. 
We represent $\GL_N((t))$ as a semidirect product 
$$\GL_N((t))=\GL_1((t))\ltimes\SL_N((t))$$ 
where $\GL_1((t))$ consists of matrices 
$\roman{diag}(a,1,...1)$. 
As for any semidirect product (see \cite{BD, Sect.~1.7}), 
the category of central extensions of $\GL_N((t))$ becomes identified 
with the category of triples $(E,E',\rho)$ where : 

\vskip2mm 

\itemitem{--} 
$E$ is a central extension of $\SL_N((t))$, 

\vskip2mm 

\itemitem{--} 
$E'$ is a central extension of $\GL_1((t))$, 

\vskip2mm 

\itemitem{--} 
$\rho$ is an action of $\GL_1((t))$ on $E$ lifting the action on 
$\SL_N((t))$ by conjugation. 

\vskip2mm 

\noindent 
Consider the triple $(E,E',\rho)$ corresponding to our extension $\tilde G$. 
By the above lemmas, both $E$ and $E'$ are trivial. 
So $\rho$ is a morphism of group ind-schemes 
$$\rho : \GL_1((t))\times\SL_N((t))\to\CC^\times.$$ 
By Lemma 5.4.5 the map $\rho$ does not depend on the second variable. 
So it is uniquely determined by its restriction to 
$\GL_1((t))\times T_\SL((t))$. This restriction, however, classifies 
the extension $\widetilde {G}|_{T((t))}$ which is trivial by Lemma 5.4.1. 
So $\rho$ is trivial and so is $\widetilde {G}$. 
Theorem 5.3.1 is proved. 

\vskip 2cm 

\centerline {\bf 6. Application to the anomaly 
of the loop space} 
\centerline {\bf and to chiral differential operators.} 

\vskip .3cm 

\subhead{(6.1) Reminder on the formal loop space}\endsubhead 
  
Let $X$ be an affine scheme over $\CC$. 
The ind-scheme $\Rc\bigl(X\otimes\CC((t))\bigr)$ will be denoted by $\widetilde{\Lc} X$. 
Thus 
$$\Hom_\Ischb(\Spec(A),\widetilde{\Lc} X)=\Hom_\Schb(\Spec(A((t))),X).\leqno(6.1.1)$$ 
Let also $\Lc^0X\subset\widetilde{\Lc} X$ be the subscheme such that 
$$\Hom_\Schb(\Spec(A),\Lc^0 X)=\Hom_\Schb(\Spec(A[[t]]),X).\leqno(6.1.2)$$ 
The formal completion of $\widetilde{\Lc} X$ along $\Lc^0X$ will be denoted by 
$\Lc X$ and called the formal loop space of $X$. 

For a commutative ring $A$ let $A((t))^\surd\subset A((t))$ 
be the subring of series $\sum_{n\in\ZZ}a_nt^n$ such that for any $n<0$ 
the element $a_n\in A$ is nilpotent. Thus 
$$\Hom_\Ischb(\Spec(A),\Lc X)=\Hom_\Schb(\Spec(A((t))^\surd),X).
\leqno(6.1.3)$$ 
It is well-known that the scheme $\Lc^0X$ can be defined for any variety $X$, 
affine or not, and the condition (6.1.2) holds. 
Further, for any scheme $S$ 
$$\Hom_\Schb(S,\Lc^0 X)=\Hom_\Lrsb((S,\Oc_S[[t]]),(X,\Oc_X)),$$ 
where $\Lrsb$ is the category of locally ringed spaces. 

It was shown in \cite{KV1} that $\Lc X$ can also be defined for an arbitrary 
$X$ 
and, in addition to (6.1.3), we have, for any scheme $S$, 
$$\Hom_\Ischb(S,\Lc X)=\Hom_\Lrsb((S,\Oc_S((t))^\surd),(X,\Oc_X)).$$ 
Thus, for any $f:S\to\Lc X$, we will denote by 
$$f_\sharp: S\to X, 
\quad 
f^\flat:f_\sharp^{-1}\Oc_X\to\Oc_S((t))^\surd$$ 
the corresponding morphism of ringed spaces. 
Thus we have a diagram 
$$X{\buildrel p\over\leftarrow}\Lc^0X{\buildrel i\over\hookrightarrow}\Lc X.$$ 
The ind-scheme $\Lc X$ being a formal neighborhood of $\Lc X$, 
it is completely determined by the sheaf of topological rings 
$\Oc_{\Lc X}$ on $(\Lc^0X)_\Zar=(\Lc X)_\Zar$. 

\subhead{(6.2) The evaluation map}\endsubhead 

Let $X$ be an arbitrary scheme over $\CC$, 
and $U\subset X$ be an affine open subset. 
We have then the evaluation map 
$$\epsilon_U:\bigl((\widetilde{\Lc} U)_\Zar,\Oc_{\tilde\Lc U}((t))\bigr)\to 
\biggl(\bigl(U\otimes\CC((t))\bigl)_\Zar,\Oc_{U\otimes\CC((t))}\biggr)$$ 
from Sect.~5. 
Consider the composition of $\epsilon_U$ with the embedding 
$$\bigl((\Lc U)_\Zar,\Oc_{\Lc U}((t))\bigr)\subset 
\bigl((\tilde\Lc U)_\Zar,\Oc_{\tilde\Lc U}((t))\bigr)$$ 
and the projection (morphism of schemes) 
$$U\otimes\CC((t))\to U.$$ 
The resulting morphism of ringed spaces will be denoted by 
$$ev_U:\bigl(\Lc^0U,\Oc_{\Lc U}((t))\bigr)\to(U,\Oc_U).$$ 
Since for $U'\subset U$ we have 
$ev_{U'}=ev_U|_{\Lc^0U'},$ 
we have that the 
$ev_U,$ $U\subset X$, 
glue together into a morphism of ringed spaces 
$$ev=ev_X:\bigl(\Lc^0X,\Oc_{\Lc X}((t))\bigr)\to(X,\Oc_X)$$ 
which we also call the evaluation map. 
Its underlying morphism of topological spaces is 
$p:\Lc^0X\to X$. 
It is the algebro-geometric analog of the evaluation map 
$$S^1\times\Map(S^1,X)\to X$$ 
for the space of smooth loops. 

\subhead{(6.3) The bundles $\Ec_\Lc$}\endsubhead 

Let $X$ be a $\CC$-scheme and $\Ec$ a quasi-coherent sheaf of $\Oc_X$-modules. 
We denote 
$$\Ec_\Lc=ev^*(\Ec)=p^{-1}\Ec\otimes_{p^{-1}\Oc_X}\Oc_{\Lc 
X}((t)).\leqno(6.3.1)$$ 
This is a $\Oc_{\Lc X}((t))$-module. 
The relation between $ev$ and $\epsilon_U$, 
$U\subset X$ affine open, allows us to describe $\Ec_\Lc$ via $\epsilon_U$. 
Namely, let $\Ec\otimes\CC((t))$ be the quasicoherent sheaf on 
$X\otimes\CC((t))$ 
obtained by extension of scalars. Then 
$$\Ec_\Lc=\epsilon_U^*\bigl(\Ec\otimes\CC((t))\bigr)|_{\Lc X}.\leqno(6.3.2)$$ 
For an ind-scheme $Z$ one defines the sheaf $\Omega^1_Z$ on $Z_\Zar$ 
as in \cite{D}. 
Then, Lemma 6.5 of $op.$ $cit.$ together with (6.3.2) imply the following 
fact. 

\proclaim{(6.3.3) Theorem} 
Let $X$ be a smooth $\CC$-variety. Then 
$\Omega^1_{\Lc X}=\Omega^1_{X,\Lc}$. 
\endproclaim 

\subhead{(6.4) The RR theorem for the bundles $\Ec_\Lc$}\endsubhead 

By definition, for every scheme $S$ and any morphism $f:S\to\Lc X$ 
the preimage of $\Ec_\Lc$ under $f$ is equal to 
$$\Ec_{\Lc,f}=f_\sharp^{-1}(\Ec)\otimes_{f_\sharp^{-1}(\Oc_X)}\Oc_S((t)).$$ 
Thus, the determinantal classes of the $\Ec_{\Lc,f}$ give rise to a class 
$$[\Det(\Ec_\Lc)]\in H^2(\Lc X,\Oc_{\Lc X}^\times).$$ 
\proclaim{(6.4.1) Theorem} 
We have that $[\Det(\Ec_\Lc)]=\partial(ev^*(2\ch_2(\Ec)))$, 
where 
$$ev^*:H^2(X,K_2(\Oc_X))\to H^2\biggl( \Lc X,K_2 \bigl( \Oc_{\Lc X}((t))\bigr) \biggr)$$ 
is the morphism induced by $ev$ on $K_2$ and 
$$\partial: 
H^2\biggl( \Lc X,K_2\bigl( \Oc_{\Lc X}((t))\bigr) \biggr)\to H^2(\Lc X,\Oc_{\Lc X}^\times)$$ 
is the boundary map. 
\endproclaim 

\noindent{\sl Proof:} 
Follows from the local Riemann-Roch (5.3.1) applied to any $(S,f)$. 
\qed 

\subhead{(6.5) Consequences for the gerbe of chiral differential operators} 
\endsubhead 

Let $X$ be a smooth $\CC$-variety. 
The category $\Db_{\Lc X}$ of right $\Dc$-modules on $\Lc X$ 
was constructed in \cite{KV1} as an abstract category, but not as a category 
of sheaves. In the forthcoming paper [KV2] it is shown that a trivialization 
$\tau$ of the gerbe $\Det(\Omega^1_{\Lc X})$ gives a functor 
$$\gamma_\tau:\Db_{\Lc X}\to\Shb_X.$$ 
In particular, applying $\gamma_\tau$ to the object 
$$i_!p^!(\omega_X)\in\Db_{\Lc X}$$ 
one gets a sheaf of chiral differential operators. 
Thus, Theorem (6.4.1) provides a conceptual explanation of the result 
of \cite{GMS1,2} relating the class of the gerbe of $CDO$ with 
$\ch_2(\Omega^1_X)$.

\vfill\eject

\vskip 2cm

\noindent M.K.: Department of Mathematics, Yale University,
10 Hillhouse Avenue, New Haven CT 06520 USA, email:
$<$mikhail.kapranov\@yale.edu$>$.

\vskip .2cm

\noindent E.V.: D\'epartement de Math\'ematiques, Universit\'e
de Cergy-Pontoise, 2 Av. A. Chauvin, 95302 Cergy-Pontoise Cedex, France,
email: $<$ eric.vasserot\@math.u-cergy.fr$>$

\bye